\documentclass[10pt]{article}
\makeatletter
\renewcommand{\thefootnote}{}

\usepackage{pdfsync}
\usepackage{cite}
\usepackage{amscd}
\usepackage{amsmath}
\usepackage{latexsym}
\usepackage{amsfonts}
\usepackage{mathtools}
\usepackage{amssymb}
\usepackage{amsthm}
\usepackage{graphicx}
\usepackage{indentfirst}
\usepackage{enumerate}
\usepackage{amstext}
\usepackage[dvipdfm,
            pdfstartview=FitH,
            CJKbookmarks=true,
            bookmarksnumbered=true,
            bookmarksopen=true,
            colorlinks=true, 
            pdfborder=001,   
            citecolor=magenta,
            linkcolor=blue,
            ]{hyperref}       

\allowdisplaybreaks[4]

\newtheorem{thm}{\textbf Theorem}[section]

\newtheorem{rem}{\textbf Remark}[section]

\newtheorem{prop}{\textbf Proposition}[section]
\newtheorem{definition}{\textbf Definition}[section]

\numberwithin{equation}{section}

\newcommand{\be}{\begin{eqnarray}}
\newcommand{\ee}{\end{eqnarray}}

\newcommand{\bes}{\begin{eqnarray*}}
\newcommand{\ees}{\end{eqnarray*}}

\def\no{\nonumber}

\def\v{\vspace{0.1in}}

\def\be{\begin{equation}}
\def\ee{\end{equation}}
\def\bee{\begin{eqnarray}}
\def\ene{\end{eqnarray}}
\def\bes{\begin{subequations}}
\def\ees{\end{subequations}}

\def\d{\displaystyle}

\renewcommand{\theequation}{\arabic{section}.\arabic{equation}}
\begin{document}

\baselineskip=13pt
\renewcommand {\thefootnote}{\dag}
\renewcommand {\thefootnote}{\ddag}
\renewcommand {\thefootnote}{ }

\pagestyle{plain}

\begin{center}
\baselineskip=16pt \leftline{} \vspace{-.3in} {\Large \bf
Orbital stability of peakon solutions for a generalized higher-order Camassa-Holm equation
} \\[0.2in]
\end{center}

\begin{center}
{\bf Guoquan Qin$^{\rm a,b}$, Zhenya Yan$^{\rm b,c,*}$, Boling Guo$^{\rm d}$}\footnote{$^{*}$Corresponding author. {\it Email address}: zyyan@mmrc.iss.ac.cn (Z. Yan). }  \\[0.1in]
{\it\small  $^a$National Center for Mathematics and Interdisciplinary Sciences, Academy of Mathematics and Systems Science, Chinese Academy of Sciences, Beijing 100190, China \\
$^b$Key Laboratory of Mathematics Mechanization, Academy of Mathematics and Systems Science, \\ Chinese Academy of Sciences, Beijing 100190, China \\
 $^c$School of Mathematical Sciences, University of Chinese Academy of Sciences, Beijing 100049, China\\
 $^d$Institute of Applied Physics and Computational Mathematics, Beijing 100088, PR China} \\
\end{center}

\vspace{0.4in}

{\baselineskip=13pt


\vspace{-0.28in}

\noindent {\bf Abstract.}\, {\small In this paper, we investigate the orbital stability issue of a generalized higher-order
Camassa-Holm (HOCH) equation, which is an higher-order extension of the quadratic CH equation.
Firstly, we show that the HOCH equation admits a global weak peakon solution by paring it with ssome smooth test function. Secondly,
with the help of two conserved quantities and the non-sgn-changing condition, we prove the orbital stability of this peakon solution in the energy space in the sense that  its shape  remains approximately the same for all times.
Our results enrich the research of the orbital stability for the CH-type equations and are useful to better understand 
the impact of higher-order nonlinearities on the dispersion dynamics.
 }

\vspace{0.15in}


 \vspace{0.1in}


\baselineskip=13pt

\section{Introduction}
\setcounter{equation}{0}

In this paper, we would like to focus on the orbital stability of the peakon solutions to the following generalized high-order
Camassa-Holm (HOCH) equation~\cite{AncoRecio2019JPAMT}
\begin{eqnarray}\label{gCH}
m_{t}+\left(u^{2}-u_{x}^{2}\right)^{n-1}u_xm
+\partial_x\left[\left(u^{2}-u_{x}^{2}\right)^{n-1} um\right]=0,\quad m=u-u_{xx},
\end{eqnarray}
where  $n\in \mathbb{N}^{+}.$
Eq.~(\ref{gCH}) was showed to admit the single peakon solution
and two conserved quantities~\cite{AncoRecio2019JPAMT}
\begin{equation}\label{Eu}
H_1(u)=\int_{\mathbb{R}}\left(u^{2}+u_{x}^{2}\right) \mathrm{d} x,\qquad
\hat H_2(u)=\int_{-\infty}^{\infty} \frac{1}{2 n}\left(u^{2}-u_{x}^{2}\right)^{n} m \mathrm{d} x.
\end{equation}
Expanding the integrand and dropping the factor $\frac{1}{2n}$ in $\hat H_2(u)$ give its equivalence as
\begin{equation}\label{Fu}
H_2(u)=\int_{\mathbb{R}}\left(u^{2 n+1}
+\sum_{k=1}^{n} \frac{(-1)^{k+1}}{2 k-1} C_{n}^{k} u^{2 n-2 k+1} u_{x}^{2 k} \right)
\mathrm{d} x.
\end{equation}
 Qu--Fu~\cite{QuFu2020JDDE} first established the local well-posedness of the corresponding Cauchy problem of Eq.~(\ref{gCH})
 in the setting of Besov spaces $B^{s}_{p,r}$ with $s>\max\{2+1/p,5/2\}$.
Chen--Deng--Qiao~\cite{ChenDengQiao2021MM}
first proved   that Eq.~(\ref{gCH}) with the special case $n=2$  admits  global
peakon and periodic peakon solutions, and then established their orbital stability
in the energy space.

At $n=1,$ Eq.~(\ref{gCH}) reduces to the remarkable CH equation~\cite{CamassaHolm1993PRL}
\begin{eqnarray}\label{CH}
m_{t}+2u_{x}m+um_{x}=0, \quad m=u-u_{xx}.
\end{eqnarray}
It is widely known that the CH equation
was originally put forward as a bi-Hamiltonian system
by  Fokas--Fuchssteiner~\cite{FokasFuchssteiner1981PD} and
 was later re-obtained  as an
approximation to the Euler equations of hydrodynamics by  Camassa-Holm~\cite{CamassaHolm1993PRL}.
It can represent  the propagation of axially symmetric waves in
hyperelastic rods~\cite{Dai1998AM,ConstantinStrauss2000PLA}.
The CH equation is completely integrable such that it admits a Lax pair, bi-Hamiltonian structure, and
infinitely many conservation laws
~\cite{CamassaHolm1993PRL,
FokasFuchssteiner1981PD,
FisherSchiff1999PLA}.
The Lax pair enables  one   to  solve the Cauchy problem of (\ref{CH})
by the inverse scattering transform (IST)~\cite{ConstantinGerdjikovIvanov2006IP}.
The local or global  well/ill-posedness of the Cauchy problem of Eq.~(\ref{CH}) has been
discussed in~\cite{Byers2006IUMJ,
ConstantinEscher1998CPAM,
Danchin2001DIE,
Constantin2000AIF,
ConstantinEscher1998AM,
LiOlver2000JDE,GuoLiuMolinetYin2019JDE}.
An important feature of the CH equation is that
it can model the wave-breaking
phenomena, namely, the solution itself remains bounded but its slope becomes unbounded in finite time~\cite{Constantin2000AIF,
ConstantinEscher1998ASNSP,
ConstantinEscher1998CPAM,
ConstantinEscher1998AM}.
This phenomena cannot be described by the classical integrable systems, such as the KdV
equation  and the Schr\"{o}dinger equation.
Another sgnificant property of the CH equation is
 that it
admits  exact weak peakon solutions of the form
$u(t, x)=c e^{-|x-c t|}, \, c \in \mathbb{R},$
which are peaked traveling waves with a discontinuous derivative at
the crest.
Peakons are a kind of solitons.
When interacting with another one, the peakon will keep
its shape and speed unchanged,
so  it is reasonable to expect the stability of the peakons,
an elementary question for nonlinear
wave equations.
Here, for CH-type equations,
the appropriate concept of stability is orbital
stability.
Namely, if a wave is  close to a solitary wave initially,
then it will remain close to some
translate of the solitary wave  at all later times,
or to put it briefly,
 the shape of the wave remains approximately
the same for all times~\cite{ConstantinStrauss2000CPAM}.
The orbital stability issue of single peakon for  Eq.~(\ref{CH})
has been tackled by Constantin--Strauss~\cite{ConstantinStrauss2000CPAM}
with the aid of the conserved densities and  specific structure of the peakons.
Constantin--Molient~\cite{ConstantinMolinet2001PD}
 introduced a variational
approach to settle  the orbital stability issue.


The nonlinearity in Eq.~(\ref{CH}) is quadratic.
However, there do exist other CH-type equations with
cubic or higher-order nonlinearities.
For example, the following cubic modified Camassa-Holm equation
\begin{equation}\label{mCH}
m_{t}+\left(\left(u^{2}-u_{x}^{2}\right) m\right)_{x}=0, \quad m=u-u_{x x},
\end{equation}
which is also called the FORQ equation~\cite{HimonasMantzavinos2014NA,YangLiZhao2018AA}.
It was first deduced by Fuchssteiner~\cite{Fuchssteiner1996PD}
and Olver-Rosenau~\cite{OlverRosenau1996PRE}
 as a new generalization of integrable system by applying
 tri-Hamiltonian duality  to the bi-Hamiltonian representation
 of the modified KdV equation.
 Later, Qiao \cite{Qiao2006JMP} discussed its
 integrability and the structure of solutions.
Like the CH equation,
it is also completely integrable
and admits a Lax pair~\cite{OlverRosenau1996PRE,QiaoLi2011TMP}.
The Cauchy problem in the Besov setting as well as the blow-up data
was established by Fu--Gui--Qu--Liu~\cite{FuGuiQuLiu2013JDE}.
 Gui--Liu--Olver--Qu~\cite{GuiLiuOlverQu2013CMP}
addressed the blow-up criterion and another type of wave-breaking data
as well as the explicit form of the single peakon and the dynamical systems
satisfied by the multi-peakons.
The form of the single peakon is
$u(t, x)=\sqrt{3c/2} e^{-|x-c t|}, \, c\in\mathbb{R}$
whose orbital stability has been tackled by Qu--Liu--Liu~\cite{QuLiuLiu2013CMP}.

An instance of  the modified CH equation with higher-order nonlinearities is written as
\begin{equation}\label{guoCH}
m_{t}+\left(\left(u^{2}-u_{x}^{2}\right)^{n} m\right)_{x}=0, \quad m=u-u_{x x},
\end{equation}
which was deduced by Anco--Recio~\cite{AncoRecio2019JPAMT}.
Its local well-posedness of the Cauchy problem in Besov spaces
was established by Yang--Li--Zhao~\cite{YangLiZhao2018AA}.
The  single peakon solution
$$
\varphi_{c}(t, x)=a e^{-|x-c t|}, \quad c=\left(1-\sum_{k=1}^{n} \frac{(-1)^{k+1}}{2 k+1} C_{n}^{k}\right) a^{2 n},
$$
was first obtained by Anco--Recio~\cite{AncoRecio2019JPAMT}, and later re-derived by
Guo--Liu--Liu--Qu~\cite{GuoLiuLiuQu2019JDE}.
By combining
 the  idea of Constantin-Strauss~\cite{ConstantinStrauss2000CPAM}
 and Qu--Liu--Liu~\cite{QuLiuLiu2013CMP},
 the authors~\cite{GuoLiuLiuQu2019JDE} also discussed
 the orbital stability of the single peakon,
where the following inequality
 \begin{equation}\label{guo1}
\frac{n\left(2-c_{1}\right)}{n+1} M^{2 n+2}-\frac{2-c_{1}}{2} M^{2 n} E(u)+\tilde F(u) \leq 0
\end{equation}
 plays a crucial role in the proof, where $M(t)=\max_{x\in \mathbb{R}}\{u(x,t)\}$, $E(u)=H_1(u)$ is given by Eq.~(\ref{Eu}), and
 \begin{equation}\label{Fu}
\tilde F(u)=\int_{\mathbb{R}}\left(u^{2(n+1)}
+\sum_{k=1}^{n} \frac{(-1)^{k+1}}{2 k-1} C_{n+1}^{k} u^{2(n-k+1)} u_{x}^{2 k}+\frac{(-1)^{n}}{2n+1}u_x^{2(n+1)} \right)
\mathrm{d} x.
\end{equation}
 Another  example  of  the CH-type equation with higher-order nonlinearities
 is just Eq.~(\ref{gCH}). Some properties  in the setting of Besov spaces $B^{s}_{p,r}$ with $s>\max\{2+1/p,5/2\}$ of Eq.~(\ref{gCH}) are listed in  {\bf Appendix A}.

Motivated by the works~\cite{ConstantinStrauss2000CPAM,GuoLiuLiuQu2019JDE,ChenDengQiao2021MM},
we will explore  the orbital stability of the single peakon of Eq.~(\ref{gCH})
in this paper. For this purpose, we invoke similar method
analogous to that employed in~\cite{ConstantinStrauss2000CPAM,GuoLiuLiuQu2019JDE}.
Of course, the two conserved quantities $E(u)$ and $F(u)$ will play a vital role  in
this issue.
 We will
prove  the following  significant inequality (see Proposition~\ref{lemGuo32})
\begin{eqnarray}\label{guo33-intro}
\frac{(2n-1)\left(2-c_{1}\right)}{2n+1} M^{2 n+1}(t)-\frac{2-c_{1}}{2} M^{2 n-1}(t) H_1(u)+H_2(u) \leq 0,
\end{eqnarray}
which
relates $H_{1,2}(u)$  and the maximal value of
approximate solutions $u$.
One can compare (\ref{guo33-intro})
with (\ref{guo1}).
To obtain (\ref{guo33-intro}),
we  first define the same functional
$g$ as that in~\cite{ConstantinStrauss2000CPAM} or  in the CH and mCH equations
due to the same conserved quantity $H_1(u)$.
Note that the constructed function $g$ needs to vanish at the peakons.
We also need to find the other function $h$ as previous literatures did.
This can be done by first assuming  the expression of $h$ with some coefficients
$c_{k}, d_{k}$ to be determined and then comparing the expansion of $\int g^{2}h$ with the formula $H_2(u)$ and finally obtaining the expression of $c_{k}, d_{k}$
by solving recursive formulas of the two sequence.
The function $h$ is also required to vanish at the peakons.
The third ingredient in establishing (\ref{guo33-intro})
is to prove
\begin{eqnarray}\label{guo38-intro}
h(t, x) \leq \frac{2-c_{1}}{2} M^{2 n-1}(t, x) \quad \text { for } \forall(t, x) \in[0, T) \times \mathbb{R}
\end{eqnarray}
with $M(t)=\max_{x\in \mathbb{R}}\{u(x,t)\}.$
Notice that in the case of the CH equation, the function $h$
satisfy $h=u\leq M$, while in the case of mCH equation, the
function $h$ satisfy $h\leq \frac{4}{3}M^{2}$.
This results from the
 higher-order conserved quantity $F(u)$ in our case.
One can also compare (\ref{guo38-intro}) with that in~\cite{GuoLiuLiuQu2019JDE},
where $h$ satisfy $h\leq \frac{2-c_{1}}{2}M^{2n}.$
This is due to the fact that the conserved quantity $H$
here is a degree lower than that in~\cite{GuoLiuLiuQu2019JDE},
where the corresponding  $\tilde{H}(u)=\frac{1}{2(n+1)}\int_{-\infty}^{\infty}  \left(u^{2}-u_{x}^{2}\right)^{n} m u\mathrm{d} x.$
To obtain (\ref{guo38-intro}), we need to carefully analyze
the properties of the coefficients $c_{k}, d_{k}$ of the function $h$
under the non-sgn-changing condition.

Having established (\ref{guo33-intro}), it is then ready to settle
the orbital stability issue.
One first expands the conserved quantity $E(u)$ around the peakon $\varphi_{c}$,
then uses the error term $\left|M-\max _{x \in \mathbb{R}} \varphi_{c}\right|$
to control  the term $|u-\varphi_{c}(\cdot-\xi(t))|_{H^{1}}$.
Finally, the  root structure of the  polynomial inequality
$Q(y)$ defined By Eq.~(\ref{guo326}) will be analyzed to estimate the error term.

To state our main results,
we need to  define the weak solutions associated with Eq.~(\ref{gCH}).
For  this purpose, we recast it as
\begin{equation}\label{guo24}
\begin{aligned}
&u_{t}-u_{t x x}
+\sum_{k=0}^{n-1}(-1)^{k}C_{n-1}^{k}u^{2 n-2 k-2} u_{x}^{2 k+1}[(2n-2k+1)u+2(2n-2k-1)u_{x x}] \\
&\qquad +2\sum_{k=0}^{n-1}(-1)^{k}   C_{n-1}^{k}(n-k-1) u^{2 n-2 k-3} u_{x}^{2 k+1} u_{x x}^{2}
+\sum_{k=0}^{n-1}(-1)^{k+1} C_{n-1}^{k} u^{2 n-2 k-1} u_{x}^{2 k} u_{x x x}=0.
\end{aligned}
\end{equation}
Then operating  $(1-\partial_{x}^{2})^{-1}$
to (\ref{guo24}) gives rise to
\begin{equation}\label{guo25}
\begin{aligned}
&u_{t}+\sum_{k=0}^{n-1} \frac{(-1)^{k}}{2 k+1} C_{n-1}^{k} u^{2 n-2 k-1} u_{x}^{2 k+1} +\left(1-\partial_{x}^{2}\right)^{-1} \partial_{x}\left(u^{2 n}
+\sum_{k=1}^{n} \frac{(-1)^{k-1}(2 n-2 k+1)}{2n(2 k-1)} C_{n}^{k} u^{2 n-2 k} u_{x}^{2 k}\right)\\
&\qquad+\left(1-\partial_{x}^{2}\right)^{-1}\left(\sum_{k=1}^{n-1} \frac{(-1)^{k+1}}{2 k+1} C_{n-1}^{k} u^{2 n-2 k-1} u_{x}^{2 k+1}\right)=0.
\end{aligned}
\end{equation}
Let $p(x)=\frac{1}{2}e^{-|x|}(x\in\mathbb{R})$.
Then there holds $p * f=\left(1-\partial_{x}^{2}\right)^{-1} f$
for  $f \in L^{2}$.
We use this combined with (\ref{guo25}) to define weak solutions
of (\ref{gCH}) as follows.

\begin{definition}\label{defGuo21}

For the given initial data $u_{0} \in W^{1,2 n}(\mathbb{R})$, $u(t, x) \in L_{l o c}^{\infty}\left([0, T) ; W_{l o c}^{1,2 n}(\mathbb{R})\right)$ is called a weak solution to the HOCH equation (\ref{gCH}) if it satisfies
$$
\begin{aligned}
&\int_{0}^{T}\!\!\! \int_{\mathbb{R}}\left[u \phi_{t}+\frac{1}{2 n} u^{2 n} \phi_{x}+\left(\sum_{k=1}^{n-1} \frac{(-1)^{k+1}}{2 k+1} C_{n-1}^{k} u^{2 n-2 k-1} u_{x}^{2 k+1}\right) \phi\right. \\
&\qquad\qquad +p *\left(u^{2 n}
+\sum_{k=1}^{n} \frac{(-1)^{k-1}(2 n-2 k+1)}{2n(2 k-1)} C_{n}^{k} u^{2 n-2 k} u_{x}^{2 k}\right) \cdot \partial_{x} \phi \\
&\qquad\qquad \left.-p *\left(\sum_{k=1}^{n-1} \frac{(-1)^{k+1}}{2 k+1} C_{n-1}^{k} u^{2 n-2 k-1} u_{x}^{2 k+1}\right) \cdot \phi\right] \mathrm{d} x \mathrm{d} t+\int_{\mathbb{R}} u_{0}(x) \phi(0, x) \mathrm{d}x=0
\end{aligned}
$$
for any test function $\phi(t, x) \in C_{c}^{\infty}([0, T) \times \mathbb{R})$. In particular, $u(t,x)$ is called a global weak solution to Eq.~(\ref{gCH}) if it is a weak solution on $[0, T)$ for every $T>0$.
\end{definition}

Then comes our first result concerning the explicit formula of the weak peakon solution
to (\ref{gCH}):

\begin{prop}\label{thmGuo21}
For any $c\in\mathbb{R}$, the peaked function of the form
\begin{equation}\label{guo26}
  u(t,x)=\varphi_{c}(t, x)=a e^{-|x-c t|}, \qquad  c=\left(1+\frac{1}{2n}\right)\sum_{k=0}^{n}\frac{(-1)^{k}}{2k+1}C_{n}^{k}a^{2n-1}
\end{equation}
is a global weak solution to equation (\ref{gCH}) in the sense of Definition \ref{defGuo21}.
\end{prop}

\begin{rem} If $u(t,x)$ is a solution of Eq.~(\ref{gCH}), then so is $-u(-t, x)$. Therefore, without loss of generalization, one can assume that $a>0$.
\end{rem}

\begin{prop}\label{lemDiffE}
For every $u \in H^{1}(\mathbb{R})$ and $\xi \in \mathbb{R}$, one has
\begin{eqnarray}
H_1(u)-H_1\left(\varphi_{c}(\cdot-\xi)\right)=\left\|u-\varphi_{c}(\cdot-\xi)\right\|_{H^{1}}^{2}+4 a(u(\xi)-a),
\end{eqnarray}
where $a$ and $c$ are defined by $c=\sum_{k=0}^{n}\frac{(-1)^{k}}{2k+1}\frac{2n+1}{2n}C_{n}^{k}a^{2n-1}$.
\end{prop}

\begin{prop}\label{lemGuo32}
 Assume $u_{0} \in H^{s}(\mathbb{R}), s>5 / 2$, and $m_{0} \geq 0 .$ Let $u(t, x)$ be the positive solution of the Cauchy problem of the HOCH equation (1.1) with initial data $u_{0} .$ Then
\begin{eqnarray}\label{guo33}
\frac{(2n-1)\left(2-c_{1}\right)}{2n+1} M^{2 n+1}(t)-\frac{2-c_{1}}{2} M^{2 n-1}(t) H_1(u)+H_2(u) \leq 0,
\end{eqnarray}
where  $M(t) \triangleq$ $\max _{x \in \mathbb{R}}\{u(t, x)\}$, and $c_1=\frac{1}{2}+\sum_{j=1}^{n}(-1)^{j+1} \frac{2 j-3}{2(2 j-1)} C_{n}^{j}$.
\end{prop}

\begin{prop}\label{lemGuo33}
 For $u \in H^{s}(\mathbb{R}), s>\frac{5}{2}$, if $\left\|u-\varphi_{c}\right\|_{H^{1}(\mathbb{R})}<\varepsilon$, with $0<\varepsilon<(\gamma-2 \sqrt{2}) a,\, \gamma>2\sqrt{2}$, then
\begin{equation}\label{ef}
\left|H_1(u)-H_1\left(\varphi_{c}\right)\right| \leq a \gamma \varepsilon, \qquad
\left|H_2(u)-H_2\left(\varphi_{c}\right)\right| \lesssim G\left(n, c,\|u\|_{H^{s}}\right)\varepsilon,
\end{equation}
where the constant $G\left(n, c,\|u\|_{H^{s}}\right)>0$ depends only on the wave speed $c$, $n\in\mathbb{N}^+$  and $\|u\|_{H^{s}}$.
\end{prop}

\begin{prop}\label{lemGuo34}
 For $0<u \in H^{s}(\mathbb{R}), s>\frac{5}{2}$, let $M=\max _{x \in \mathbb{R}}\{u(x)\}.$ If $H_{1,2}(u)$ satisfy Eq.~(\ref{ef})
then
$$
|M-a| \lesssim \sqrt{G(n, c,\|u\|_{H^{s}})\varepsilon}.
$$
\end{prop}

We then establish the orbital stability of the peakon solutions (\ref{guo26}),
that is
\begin{thm}\label{thmGuo31}
 The peakon solution $\varphi_{c}(t, x)$ defined in (\ref{guo26}) with the travelling wave speed $c$ is orbitally stable in the following sense. If $u_{0}(x) \in H^{s}(\mathbb{R})$, for some $s>5 / 2, m_{0}(x)=\left(1-\partial_{x}^{2}\right) u_{0} \neq 0$ is nonnegative, and
\begin{eqnarray*}
\left\|u(0, \cdot)-\varphi_{c}\right\|_{H^{1}(\mathbb{R})}<\varepsilon, \quad \text { for } \quad 0<\varepsilon<(\gamma-2 \sqrt{2}) a, \quad \gamma>2\sqrt{2},\,\, a>0
\end{eqnarray*}
 Then the corresponding solution $u(t, x)$ of equation (\ref{gCH}) satisfies
\begin{eqnarray*}
\sup _{t \in[0, T)}\left\|u(t, \cdot)-\varphi_{c}(\cdot-\xi(t))\right\|_{H^{1}(\mathbb{R})} \lesssim \sqrt{ a \gamma \varepsilon+4 a \sqrt{G\left(c,\left\|u_{0}\right\|_{H^{s}}\right) \varepsilon}},
\end{eqnarray*}
where $T>0$ is the maximal existence time, $\xi(t) \in \mathbb{R}$ is the maximum point of function $u(t, \cdot)$ and the constant $G\left(n, c,\left\|u_{0}\right\|_{H^{s}}\right)>0$ depends only on the wave speed $c$, $n\in\mathbb{N}^+$ and the norm $\left\|u_{0}\right\|_{H^{s}}$.
\end{thm}

The rest of this paper is organized as follows.
In section \ref{sectionPeakon}, we verify that
(\ref{guo26}) is a global weak solution to equation (\ref{gCH}) in the sense of Definition \ref{defGuo21} and show $F(u)$ is independent of time. Section \ref{sectionStability} will be used to prove the
orbital stability result--Theorem \ref{thmGuo31} by showing Propositions~\ref{lemDiffE}-\ref{lemGuo34}.

\textbf{Notation.}  $A \lesssim B$ represents $A \leq C B$ for some constant $C>0$.

\section{Peakons and conserved quantities}\label{sectionPeakon}

In this section we give the proof of Proposition \ref{thmGuo21}
After that, we will prove the quantity $F(u)$ is independent of time.

\begin{proof}[\textbf{Proof of Proposition \ref{thmGuo21}:}] We consider the peakon solution of Eq.~(\ref{gCH}) in the form
\bee
u(t,x)=\varphi_{c}(t, x)=a e^{-|x-c t|},
\ene
where $a,\, c\in\mathbb{R}$ are constants to be determined later.

First, in light of  integration by parts,
one easily verifies
\begin{equation}\label{guo27}
\partial_{x} \varphi_{c}(t, x)=-\operatorname{sgn}(x-c t) \varphi_{c}(t, x),\qquad
\partial_{t} \varphi_{c}(t, x)=c\, \partial_{x} \varphi_{c}(t, x) \in L^{\infty}
\end{equation}
for all $t\geq0$ in the sense of distribution.
Set $\varphi_{0, c}(x) \triangleq \varphi_{c}(0, x)$.
 Then there holds
\begin{equation}\label{guo28}
\lim _{t \rightarrow 0^{+}}\left\|\varphi_{c}(t, \cdot)-\varphi_{0, c}(x)\right\|_{W^{1, \infty}}=0.
\end{equation}
Invoking  (\ref{guo27})-(\ref{guo28}),
one derives   for any test function $\phi(t, x) \in C_{c}^{\infty}([0, \infty) \times \mathbb{R})$  that
\begin{equation}\label{guo210}
\begin{aligned}
&\int_{0}^{+\infty}\!\!\!\! \int_{\mathbb{R}}\left[ \varphi_{c}\partial_{t} \phi
+\frac{1}{2n} \varphi_{c}^{2n}\phi_{x}
+\left(\sum_{k=1}^{n-1} \frac{(-1)^{k+1}}{2 k+1} C_{n-1}^{k} \varphi_{c}^{2 n-2 k-1}\left(\partial_{x} \varphi_{c}\right)^{2 k+1}\right) \phi\right] \mathrm{d} x \mathrm{d} t +\int_{\mathbb{R}} \varphi_{0, c}(x) \phi(0, x) d x \\
&\quad=-\int_{0}^{+\infty}\!\!\!\! \int_{\mathbb{R}} \left[\partial_{t}\varphi_{c}
+ \varphi_{c}^{2n-1}\partial_{x}\varphi_{c}
-\left(\sum_{k=1}^{n-1} \frac{(-1)^{k+1}}{2 k+1} C_{n-1}^{k} \varphi_{c}^{2 n-2 k-1}\left(\partial_{x} \varphi_{c}\right)^{2 k+1}\right) \right]\phi \mathrm{d} x \mathrm{d} t\\
&\quad=-\int_{0}^{+\infty}\!\!\!\! \int_{\mathbb{R}}\phi \operatorname{sgn}(x-c t) \varphi_{c}\left[c-\left(1-\sum_{k=1}^{n-1} \frac{(-1)^{k+1}}{2 k+1} C_{n-1}^{k}\right) \varphi_{c}^{2 n-1}\right]  \mathrm{d} x \mathrm{d} t.
\end{aligned}
\end{equation}
Substituting (\ref{guo26}) into the above integrands,
one derives that
\begin{equation}\label{guo211}
\begin{aligned}
\operatorname{sgn}(x\!-\!c t) \varphi_{c}\!\left[\!c\!-\!\left(1-\sum_{k=1}^{n-1} \frac{(-1)^{k+1}}{2 k+1} C_{n-1}^{k}\right) \varphi_{c}^{2 n-1}\!\right]\! =\!a^{2 n}\left(1\!-\!\sum_{k=1}^{n-1} \frac{(-1)^{k+1}}{2 k+1} C_{n-1}^{k}\right)\left(e^{c t-x}
\!-\!e^{2 n(c t\!-\!x)}\right)
\end{aligned}
\end{equation}
for  $x>c t$ and
\begin{equation}\label{guo212}
\begin{aligned}
\operatorname{sgn}(x\!-\!c t) \varphi_{c}\!\left[\!c\!-\!\left(1\!-\!\sum_{k=1}^{n-1} \frac{(-1)^{k+1}}{2 k+1} C_{n-1}^{k}\!\right) \varphi_{c}^{2 n-1}\!\right]\!=\!-\!a^{2 n}\left(1\!-\!\sum_{k=1}^{n-1} \frac{(-1)^{k+1}}{2 k+1} C_{n-1}^{k}\!\right)\!
\left(e^{x\!-\!c t}\!-\!e^{2 n(x\!-\!c t)}\right)
\end{aligned}
\end{equation}
for  $x\leq c t$, where one can easily verify that
$\sum_{k=0}^{n}\frac{(-1)^{k}}{2k+1}\frac{2n+1}{2n}C_{n}^{k}
=1-\sum_{k=1}^{n-1} \frac{(-1)^{k+1}}{2 k+1} C_{n-1}^{k}$.

Substituting $\varphi_{c}$ in the remainder terms in Definition \ref{defGuo21}
and employing (\ref{guo27}) lead to
\begin{eqnarray}\label{guo213}
&&\int_{0}^{+\infty}\!\!\!\! \int_{\mathbb{R}}\bigg[
p *\left(\varphi_{c}^{2 n}
+\sum_{k=1}^{n} \frac{(-1)^{k-1}(2 n-2 k+1)}{2n(2 k-1)} C_{n}^{k} \varphi_{c}^{2 n-2 k} (\partial_{x}\varphi_{c})^{2 k}\right) \cdot \partial_{x}\phi \nonumber\\
&& \quad\qquad -p *\left(\sum_{k=1}^{n-1} \frac{(-1)^{k+1}}{2 k+1} C_{n-1}^{k} \varphi_{c}^{2 n-2 k-1} (\partial_{x}\varphi_{c})^{2 k+1}\right) \cdot \phi\bigg] \mathrm{d} x \mathrm{d} t\nonumber\\
&&=-\int_{0}^{+\infty}\!\!\!\! \int_{\mathbb{R}}\bigg[
p_{x} *\left(
\sum_{k=1}^{n} \frac{(-1)^{k-1}(2 n-2 k+1)}{2n(2 k-1)} C_{n}^{k} \varphi_{c}^{2 n-2 k} (\partial_{x}\varphi_{c})^{2 k}\right) \cdot \phi \nonumber\\
&&\quad\qquad +p *\left(2n\varphi_{c}^{2 n-1}\partial_{x}\varphi_{c}+\sum_{k=1}^{n-1} \frac{(-1)^{k+1}}{2 k+1} C_{n-1}^{k} \varphi_{c}^{2 n-2 k-1} (\partial_{x}\varphi_{c})^{2 k+1}\right) \cdot \phi\bigg] \mathrm{d} x \mathrm{d} t \no\\
&&=-\int_{0}^{+\infty}\!\!\!\! \int_{\mathbb{R}}
 \phi \cdot p_{x} *\bigg[
\sum_{k=1}^{n} \frac{(-1)^{k-1}(2 n-2 k+1)}{2n(2 k-1)} C_{n}^{k} \varphi_{c}^{2 n-2 k} (\partial_{x}\varphi_{c})^{2 k}   +\left(1+\sum_{k=1}^{n-1} \frac{(-1)^{k+1}}{2n(2k+1)} C_{n-1}^{k}\right) \varphi_{c}^{2 n}
\bigg] \mathrm{d} x \mathrm{d} t. \qquad \no
\end{eqnarray}
Using $\partial_{x} p(x)=-\frac12\operatorname{sgn}(x) e^{-|x|}$ for $x \in \mathbb{R}$,
one finds the integrand in the above(after dropping the factor  $\phi$) is
\begin{eqnarray}\label{guo215}
-\frac{a^{2n}}{2}\bigg(1+\sum_{k=1}^{n-1}\frac{(-1)^{k+1}}{2n(2k+1)} C_{n-1}^{k}\!+\!\sum_{k=1}^{n} \frac{(-1)^{k-1}(2 n-2 k+1)}{2n(2 k-1)} C_{n}^{k}\bigg)\int_{\mathbb{R}}\operatorname{sgn}(x\!-\!y)e^{-|x-y|-2 n|y-c t|} \mathrm{d} y
\end{eqnarray}
When $x>ct$, we have
\begin{eqnarray}\label{guo216}
&&(\ref{guo215})
=-\frac{a^{2n}}{2}\bigg(1+\sum_{k=1}^{n-1}\frac{(-1)^{k+1}}{2n(2k+1)} C_{n-1}^{k}
+\sum_{k=1}^{n} \frac{(-1)^{k-1}(2 n-2 k+1)}{2n(2 k-1)} C_{n}^{k}\bigg)\nonumber\\
&&\qquad\qquad\qquad\quad \times
\bigg(\int_{-\infty}^{ct}+\int_{ct}^{x}+\int_{x}^{+\infty}\bigg)\operatorname{sgn}(x-y)
e^{-|x-y|-2 n|y-c t|} \mathrm{d} y.
\end{eqnarray}
Direct computation yields
\begin{eqnarray*}
\begin{array}{l}
\d \int_{-\infty}^{ct}\operatorname{sgn}(x-y)e^{-|x-y|-2 n|y-c t|} \mathrm{d} y
=\frac{1}{2n+1}e^{ct-x}, \v\\
\d\int_{ct}^{x}\operatorname{sgn}(x-y)e^{-|x-y|-2 n|y-c t|} \mathrm{d} y
=\frac{1}{2n-1}[e^{ct-x}-e^{2n(ct-x)}], \v\\
\d \int_{x}^{+\infty}\operatorname{sgn}(x-y)e^{-|x-y|-2 n|y-c t|} \mathrm{d} y
=-\frac{1}{2n+1}e^{2n(ct-x)}.
\end{array}
\end{eqnarray*}
Plugging the above into (\ref{guo216}) yields for $\pm (x-ct)\geq 0$
\begin{eqnarray}\label{guo217}
(\ref{guo215})=\frac{a^{2n}}{\pm(1-4n^2)}\bigg( 2n+\sum_{k=1}^{n-1}
\frac{(-1)^{k+1}}{2 k+1} C_{n-1}^{k}
+\sum_{k=1}^{n} \frac{(-1)^{k-1}(2 n-2 k+1)}{2 k-1} C_{n}^{k}\bigg)[e^{\pm(ct-x)}-e^{\pm 2n(ct-x)}].
\end{eqnarray}

Moreover one can prove the identity
\begin{eqnarray}\label{guo220}
\frac{1}{4n^2-1}\left(2 n+\sum_{k=1}^{n-1}
\frac{(-1)^{k+1}}{2 k+1} C_{n-1}^{k}
+\sum_{k=1}^{n} \frac{(-1)^{k-1}(2 n-2 k+1)}{2k-1} C_{n}^{k}\right)=1-\sum_{k=1}^{n-1} \frac{(-1)^{k+1}}{2 k+1} C_{n-1}^{k}
\end{eqnarray}
according to the relations~\cite{GuoLiuLiuQu2019JDE}
\begin{eqnarray}\no
\sum_{k=1}^{n} \frac{(-1)^{k-1}}{2 k-1} C_{n}^{k}
=\frac{(2n)!!}{(2n-1)!!}-1,\quad 
\sum_{k=0}^{n} \frac{(-1)^{k}}{2 k+1} C_{n}^{k}
=1-\sum_{k=1}^{n} \frac{(-1)^{k+1}}{2 k+1} C_{n}^{k}
=\frac{(2n)!!}{(2n+1)!!}.
\end{eqnarray}

Gathering (\ref{guo211})--(\ref{guo212}) and (\ref{guo217})--(\ref{guo220}),
one deduces for any $\phi(t, x) \in C_{c}^{\infty}([0, \infty) \times \mathbb{R})$ that
\begin{equation*}
\begin{aligned}
&\int_{0}^{+\infty}\!\!\!\! \int_{\mathbb{R}}\left[\varphi_{c} \phi_{t}+\frac{1}{2 n} \varphi_{c}^{2 n} \phi_{x}+\left(\sum_{k=1}^{n-1} \frac{(-1)^{k+1}}{2 k+1} C_{n-1}^{k} \varphi_{c}^{2 n-2 k-1} (\partial_{x}\varphi_{c})^{2 k+1}\right) \phi\right. \\
&\qquad +p *\left(\varphi_{c}^{2 n}
+\sum_{k=1}^{n} \frac{(-1)^{k-1}(2 n-2 k+1)}{2n(2 k-1)} C_{n}^{k} \varphi_{c}^{2 n-2 k} (\partial_{x}\varphi_{c})^{2 k}\right) \cdot \partial_{x} \phi \\
&\qquad\left.-p *\left(\sum_{k=1}^{n-1} \frac{(-1)^{k+1}}{2 k+1} C_{n-1}^{k} \varphi_{c}^{2 n-2 k-1} (\partial_{x}\varphi_{c})^{2 k+1}\right) \cdot \phi\right] \mathrm{d} x \mathrm{d} t+\int_{\mathbb{R}} \varphi_{0,c}(x) \phi(0, x) \mathrm{d}x=0.
\end{aligned}
\end{equation*}
We thus complete the proof of  Proposition \ref{thmGuo21}.
\end{proof}

We next  prove that the quantity $F(u)$ given by Eq.~(\ref{Fu}) is independent of time.

To do this, we
set $v(x, t)=\int_{-\infty}^{x} u_{t}(z, t) \mathrm{d} z$ and employ  integration by parts to obtain
\begin{equation}\label{guo47}
\begin{aligned}
&\frac{\mathrm{d}}{\mathrm{d} t} \int_{\mathbb{R}} u^{2 n+1} \mathrm{d} x
=(2 n+1) \int_{\mathbb{R}} u^{2 n} v_{x} \mathrm{d} x
=-(2 n+1) \int_{\mathbb{R}}(2 n) u^{2 n-1} u_{x} v \mathrm{d} x,
\end{aligned}
\end{equation}
and
\begin{equation}\label{guo49}
\begin{aligned}
\frac{\mathrm{d}}{\mathrm{d} t} \int_{\mathbb{R}}n u^{2 n-1} u_{x}^{2} \mathrm{d} x
=&\,n \int_{\mathbb{R}}\left[(2n-1)u^{2n-2}u_{x}^{2}v_{x}+2u^{2n-1}u_{x}v_{xx}\right ]\mathrm{d} x\\
=&\, n \int_{\mathbb{R}}\left[(2 n-1)(2 n-2) u^{2 n-3} u_{x}^{3}
+4(2 n-1) u^{2 n-2} u_{x} u_{x x}
+2u^{2 n-1} u_{x x x}\right] v \mathrm{d} x.
\end{aligned}
\end{equation}
In a similar manner,
\begin{eqnarray}\label{guo410}
&& \frac{\mathrm{d}}{\mathrm{d} t} \int_{\mathbb{R}} \sum_{k=2}^{n} \frac{(-1)^{k+1}}{2 k-1} C_{n}^{k} u^{2 n-2 k+1} u_{x}^{2 k} \mathrm{d} x \nonumber\\
&&\quad= \sum_{k=2}^{n} \frac{(-1)^{k+1}}{2 k-1} C_{n}^{k} \int_{\mathbb{R}}\left[(2n-2k+1) u^{2 n-2 k} u_{x}^{2 k} u_{t}+2 k u^{2n-2k+1} u_{x}^{2 k-1} u_{t x}\right] \mathrm{d} x \nonumber\\
&&\quad = \sum_{k=2}^{n} \frac{(-1)^{k+1}}{2 k-1} C_{n}^{k} \int_{\mathbb{R}}\left[(2n-2k+1) u^{2 n-2 k} u_{x}^{2 k} v_{x}+2 k u^{2n-2k+1} u_{x}^{2 k-1} v_{x x}\right] \mathrm{d} x \nonumber\\
&&\quad = \sum_{k=2}^{n}(-1)^{k+1} C_{n}^{k} \int_{\mathbb{R}}\Big[(2n-2k+1)(2 n-2 k) u^{2 n-2 k-1} u_{x}^{2 k+1} \no\\
&&\quad \quad +4 k(2n-2k+1)  u^{2 n-2 k} u_{x}^{2 k-1} u_{x x} +2 k(2 k-2) u^{2 n-2 k+1} u_{x}^{2 k-3} u_{x x}^{2}+2 k u^{2 n-2 k+1} u_{x}^{2 k-2}u_{xxx}\Big]v\mathrm{d} x.\quad\quad
\end{eqnarray}
Combining (\ref{guo47})-(\ref{guo410}), one finds
\begin{eqnarray}\label{Fug}
\frac{\mathrm{d}F(u)}{\mathrm{d}t}\!\!\!\!\!\!\!\!\!\!&&=\int_{\mathbb{R}}\bigg[
\sum_{k=2}^{n-1}(-1)^{k+1} C_{n}^{k}(2n-2k+1)(2n-2k)u^{2 n-2 k-1}u_{x}^{2 k+1}\nonumber\\
&&
\quad\quad+\sum_{k=1}^{n-1}(-1)^{k} C_{n}^{k+1}4(k+1)(2n-2k-1)u^{2 n-2 k-2}u_{x}^{2 k+1}u_{xx}\nonumber\\
&&
\quad\quad+\sum_{k=0}^{n-2}(-1)^{k+1} C_{n}^{k+2}2(k+2)(2k+2)u^{2 n-2 k-3}u_{x}^{2 k+1}u_{xx}^{2}\nonumber\\
&&
\quad\quad+\sum_{k=1}^{n-1}(-1)^{k} C_{n}^{k+1}2(k+1)u^{2 n-2 k-1}u_{x}^{2 k}u_{xxx}
\bigg]v\mathrm{d} x\nonumber\\&&
=\int_{\mathbb{R}}\bigg[
\sum_{k=0}^{n-1}(-1)^{k+1} C_{n-1}^{k}(2n-2k+1) 2n u^{2 n-2 k-1}u_{x}^{2 k+1}\nonumber\\&&
\quad\quad+\sum_{k=0}^{n-1}(-1)^{k} C_{n-1}^{k}4n(2n-2k-1)u^{2 n-2 k-2}u_{x}^{2 k+1}u_{xx}\nonumber\\&&
\quad\quad+\sum_{k=0}^{n-1}(-1)^{k+1} C_{n-1}^{k}4n(n-k-1)u^{2 n-2 k-3}u_{x}^{2 k+1}u_{xx}^{2}\nonumber\\&&
\quad\quad+\sum_{k=0}^{n-1}(-1)^{k} C_{n-1}^{k}2nu^{2 n-2 k-1}u_{x}^{2 k}u_{xxx}
\bigg]v\mathrm{d} x,
\end{eqnarray}
which along with (\ref{guo24}) leads to
\begin{equation*}
\frac{\mathrm{d} F(u)}{\mathrm{d} t}
=2n \int_{\mathbb{R}}\left(u_{t}-u_{t x x}\right) v \mathrm{d} x
=\int_{\mathbb{R}}2n\left(v v_{x}-v v_{x x x}\right) \mathrm{d}x=0.
\end{equation*}

\section{Orbital stability of peakon solutions}\label{sectionStability}

In this section, we prove the orbital stability of  peakons for
Eq.~(\ref{gCH}).
The proof will be divided into several lemmas.
One should first notice that
under  the assumption
of  Proposition \ref{thmGuo31},
 Eq.~(\ref{gCH}) admits a  unique local positive solution of
 by  the use of Lemmas \ref{lemGuo21}-\ref{lemGuo22}.
Also, it is easily to see that
the function
$\varphi_{c}(x)=a\varphi(x)=a e^{-|x|}\in H^{1}(\mathbb{R})$
takes its maximum at $x=0$, i.e.,
\begin{eqnarray}
\max _{x \in \mathbb{R}}\left\{\varphi_{c}(x)\right\}=\varphi_{c}(0)=a.
\end{eqnarray}
Simple calculation yields
\begin{equation}\label{Epeakon}
E\left(\varphi_{c}\right)=\left\|\varphi_{c}\right\|_{H^{1}}^{2}=a^{2} \int_{\mathbb{R}}\left(\varphi^{2}+\varphi_{x}^{2}\right) d x=2 a^{2}
\end{equation}
and
\begin{equation}\label{Fpeakon}
\begin{aligned}
F\left(\varphi_{c}\right) &=a^{2 n+1} \int_{\mathbb{R}}
\left(\varphi^{2 n+1}+\sum_{k=1}^{n} \frac{(-1)^{k+1}}{2 k-1} C_{n}^{k} \varphi^{2 n-2 k+1} \varphi_{x}^{2 k}\right) \mathrm{d} x =\frac{2a^{2 n+1}}{2n+1}\left(1+\sum_{k=1}^{n} \frac{(-1)^{k+1}}{2 k-1} C_{n}^{k}\right),
\end{aligned}
\end{equation}
where $a$ is determined implicitly by $c=\sum_{k=0}^{n}\frac{(-1)^{k}}{2k+1}\frac{2n+1}{2n}C_{n}^{k}a^{2n-1}$
provided by Theorem \ref{thmGuo21}.

\begin{proof}[\textit{\textbf{Proof of Proposition \ref{lemDiffE}}}]
First, one easily verifies   $\varphi-\partial_{x}^{2} \varphi=2 \delta$ with $\delta$ being  the Dirac distribution. Consequently, we find with the help of
integration by parts and (\ref{Epeakon}) that
$$
\begin{aligned}
\left\|u-\varphi_{c}(\cdot-\xi)\right\|_{H^{1}}^{2}=&
\int_{\mathbb{R}}\left(u^{2}\!+\!u_{x}^{2}\right) d x\!+\!\int_{\mathbb{R}}\left(\varphi_{c}^{2}\!+\!
\left(\partial_{x} \varphi_{c}\right)^{2}\right) d x\! -\! 2 a \int_{\mathbb{R}} u_{x}(x) \varphi_{x}(x-\xi) d x\! -
\!2 a \int_{\mathbb{R}} u(x) \varphi(x-\xi) d x \\
=&\, E(u)+E\left(\varphi_{c}(\cdot-\xi)\right)-2 a \int_{\mathbb{R}}\left(1-\partial_{x}^{2}\right) \varphi(x-\xi) u(x) d x \\
=&\, E(u)+E\left(\varphi_{c}(\cdot-\xi)\right)-4 a \int_{\mathbb{R}} \delta(x-\xi) u(x) d x \\
=&\, E(u)+E\left(\varphi_{c}(\cdot-\xi)\right)-4 a u(\xi)\\
=&\, E(u)-E\left(\varphi_{c}(\cdot-\xi)\right)-4 a(u(\xi)-a).
\end{aligned}
$$
This completes the proof of this Proposition \ref{lemDiffE}.
\end{proof}

\begin{proof}[\textit{\textbf{Proof of Proposition \ref{lemGuo32}}}]

Let the maximus of $u(t,x)$ be taken at $x=\xi(t)$, i.e.,
$M(t)=u(t,\xi(t))$. We  define the same function $g$
as in \cite{ConstantinStrauss2000CPAM}
\begin{eqnarray}
g(t, x) \triangleq\left\{\begin{array}{ll}
u(t, x)-u_{x}(t, x), & x<\xi(t), \v \\
u(t, x)+u_{x}(t, x), & x>\xi(t).
\end{array}\right.
\end{eqnarray}
Then it is easy to know that
\begin{eqnarray}\label{guo34}
\int_{\mathbb{R}} g^{2}(t, x) \mathrm{d} x=E(u)-2 M^{2}(t).
\end{eqnarray}

Let a function $h(t,x)$ be
\begin{equation*}
h(t, x) \triangleq\left\{\begin{array}{ll}
\left(u^{2 n-1}+\sum_{k=1}^{2 n-2} c_{k} u^{2 n-1-k} u_{x}^{k}\right)(t, x), & x<\xi(t), \v\\
\left(u^{2 n-1}+\sum_{k=1}^{2 n-2} d_{k} u^{2 n-1-k} u_{x}^{k} \right)(t, x), & x>\xi(t),
\end{array}\right.
\end{equation*}
where $c_{k}, d_{k}\, (k=1,2, \cdots, 2n-2)$, are constants given by
\begin{eqnarray}\label{guo35}
\begin{cases}
c_{1}=-d_{1}=\d\frac{1}{2}+\sum_{j=1}^{n}(-1)^{j+1} \frac{2 j-3}{2(2 j-1)} C_{n}^{j}, \v\\
c_{2 m}=d_{2 m}=\d\sum_{j=m+1}^{n}(-1)^{j+1} \frac{2 j-(2 m+1)}{2 j-1} C_{n}^{j}, \quad m=1,2, \cdots, n-1, \v\\
c_{2 m-1}=-d_{2 m-1}=\d\sum_{j=m+1}^{n}(-1)^{j+1} \frac{2 (j-m)}{2 j-1} C_{n}^{j}, \quad m=1,2 \cdots, n-1.
\end{cases}
\end{eqnarray}
With $c_{k}$ and  $d_{k}(k=1,2, \cdots, 2 n-2)$ defined above, one easily verifies that
\begin{eqnarray}\label{guo43}
\begin{cases}
c_{2 n-2}-\dfrac{(-1)^{n+1}}{2 n-1}=0, \v\\
c_{2 n-3}-2 c_{2 n-2}=0,  \v\\
 c_{k}-2 c_{k-1}+c_{k-2}=0, \quad\quad k=2 j+1,(j=1,2, \cdots, n-2), \v\\
  c_{k}-2 c_{k-1}+c_{k-2}=\dfrac{(-1)^{j+1}}{2 j-1} C_{n}^{j}, \quad k=2 j,(j=2,3, \cdots, n-1), \v\\
   c_{2}-2 c_{1}+1=C_{n}^{1}
\end{cases}
\end{eqnarray}
and
\begin{eqnarray}\label{guo44}
\begin{cases}
\dfrac{(-1)^{n+1}}{2 n-1}+d_{2 n-2}=0, \v\\
d_{2 n-3}+2 d_{2 n-2}=0,  \v\\
 d_{k}+2 d_{k-1}+d_{k-2}=0, \quad\quad k=2 j+1,(j=1,2, \cdots, n-2), \v\\
  d_{k}+2 d_{k-1}+d_{k-2}=\dfrac{(-1)^{j+1}}{2 j-1} C_{n}^{j}, \quad k=2 j,(j=2,3, \cdots, n-1), \v\\
   d_{2}+2 d_{1}+1=C_{n}^{1}.
\end{cases}
\end{eqnarray}
Simple calculation then yields
\begin{equation}\label{guo41}
\begin{aligned}
&\left(u^{2}-2 u u_{x}+u_{x}^{2}\right)\left(u^{2 n-1}
+\sum_{k=1}^{2 n-2} c_{k} u^{2 n-1-k} u_{x}^{k}\right) \\
&\quad= u^{2 n+1}+\left(c_{1}-2\right) u^{2 n} u_{x}+\left(c_{2}-2 c_{1}+1\right) u^{2 n-1} u_{x}^{2}+\sum_{k=3}^{2 n-2}\left(c_{k}-2 c_{k-1}+c_{k-2}\right) u^{2 n-k+1} u_{x}^{k} \\
&\qquad+\left(c_{2 n-3}-2c_{2 n-2}\right) u^{2} u_{x}^{2 n-1}
+c_{2 n-2}u u_{x}^{2 n}
\end{aligned}
\end{equation}
and
\begin{equation}\label{guo42}
\begin{aligned}
&\left(u^{2}+2 u u_{x}+u_{x}^{2}\right)\left(u^{2 n-1}
+\sum_{k=1}^{2 n-2} d_{k} u^{2 n-1-k} u_{x}^{k}\right) \\
&\quad= u^{2 n+1}+\left(d_{1}+2\right) u^{2 n} u_{x}+\left(d_{2}+2 d_{1}+1\right) u^{2 n-1} u_{x}^{2}+\sum_{k=3}^{2 n-2}\left(d_{k}+2 d_{k-1}+d_{k-2}\right) u^{2 n-k+1} u_{x}^{k} \\
&\qquad +\left(d_{2 n-3}+2d_{2 n-2}\right) u^{2} u_{x}^{2 n-1}
+d_{2 n-2}u u_{x}^{2 n}.
\end{aligned}
\end{equation}
Combining (\ref{guo43})-(\ref{guo42}), one finds
\begin{eqnarray}\label{guo36}
&&\int_{\mathbb{R}} h(t, x) g^{2}(t, x) \mathrm{d} x \nonumber\\
&&=\int_{-\infty}^{\xi}\!\!\left(u-u_{x}\right)^{2}\left(u^{2 n-1}\!+\!\sum_{k=1}^{2 n-2} c_{k} u^{2 n-1-k} u_{x}^{k}\right) \mathrm{d} x\!+\!
\int_{\xi}^{\infty}\!\left(u+u_{x}\right)^{2}\!\left(u^{2 n-1}\!+\!\sum_{k=1}^{2 n-2} d_{k} u^{2 n-1-k} u_{x}^{k}\right)\! \mathrm{d} x \nonumber\\
&&= \int_{-\infty}^{\xi}\left(u^{2 n+1}+\sum_{k=1}^{n} \frac{(-1)^{k+1}}{2 k-1} C_{n}^{k} u^{2 n-2 k+1} u_{x}^{2 k}\right) \mathrm{d} x \nonumber\\
&&\quad+\left(c_{1}\!-\!2\right)\! \int_{-\infty}^{\xi}\! u^{2 n} u_{x} \mathrm{d} x\!+\!\int_{\xi}^{\infty}\!\left(u^{2 n+1}\!+\!\sum_{k=1}^{n}\!\frac{(-1)^{k+1}}{2 k-1} C_{n}^{k} u^{2 n\!-\!2 k\!+\!1} u_{x}^{2 k}\right) \mathrm{d} x\!+\!\left(d_{1}\!+\!2\right) \int_{\xi}^{\infty}\!u^{2n}u_{x} \mathrm{d} x \nonumber\\
&&= F(u)+\frac{2(c_{1}-2)}{2n+1} u^{2 n+1}(t, \xi)\nonumber\\
&&=F(u)-\frac{2(2-c_{1})}{2n+1} M^{2 n+1}(t).
\end{eqnarray}

It follows from Lemma \ref{lemGuo22} that one knows
\begin{eqnarray}\label{guo37}
u(t, x) \geq 0, \text { and }\left(u \pm u_{x}\right)(t, x) \geq 0 \quad \text { for } \forall(t, x) \in[0, T) \times \mathbb{R}.
\end{eqnarray}
We will next prove
\begin{eqnarray}\label{guo38}
h(t, x) \leq \frac{2-c_{1}}{2} u^{2 n-1}(t, x) \quad \text { for } \forall(t, x) \in[0, T) \times \mathbb{R},
\end{eqnarray}
or equivalently,
\begin{eqnarray}\label{guo39}
\sum_{k=1}^{2 n-2} c_{k} u^{2 n-1-k} u_{x}^{k} \leq -\frac{c_{1}}{2} u^{2 n-1},\quad
\sum_{k=1}^{2 n-2} d_{k} u^{2 n-1-k} u_{x}^{k} \leq-\frac{c_{1}}{2} u^{2 n-1}.
\end{eqnarray}
Let  $z=u_{x}/u$. Then it suffices to show the nonpositivity of the following function $f(z)$
\begin{eqnarray}\label{guo311}
f(z)=\sum_{k=1}^{2 n-2} c_{k} z^{k}+\frac{c_{1}}{2} \leq 0, \quad \text { for } z \in[-1,1]
\end{eqnarray}
recalling (\ref{guo37}).

Using  $c_{2 j+1}-2 c_{2 j}+c_{2 j-1}=0\,(j=1,2, \cdots, n-2)$
and $c_{2 n-2}=\frac{c_{2 n-3}}{2}$  stated  in (\ref{guo43}),
we recast  $f(z)$ as
\begin{equation*}
\begin{aligned}
f(z) &=\sum_{j=1}^{n-1}\left(c_{2 j} z^{2 j}+c_{2 j-1} z^{2 j-1}\right)+\frac{c_{1}}{2} \\
&=\frac{c_{2 n-3}}{2} z^{2 n-2}
+\sum_{j=1}^{n-2}\left(\frac{c_{2 j+1}}{2} z^{2j}
+\frac{c_{2j-1}}{2} z^{2 j}\right)+\sum_{j=1}^{n-1}c_{2 j-1} z^{2 j-1}+\frac{c_{1}}{2} \\
&=\frac{(1+z)^{2}}{2} \sum_{k=1}^{n-1} c_{2 k-1} z^{2 k-2}.
\end{aligned}
\end{equation*}
Accordingly, the proof of  (\ref{guo311}) is equivalent to
\begin{eqnarray}\label{guo312}
\phi(z)=\sum_{k=1}^{n-1} c_{2 k-1} z^{2 k-2} \leq 0, \quad z \in[-1,1].
\end{eqnarray}
Since $\phi(z)$ is  even  and continuous at $z=1$, one only needs to prove  $(\ref{guo312})$  for $z \in[0,1)$.
Using the expression  of $c_{2k-1}(k=1,2, \cdots, n-1)$ provided in (\ref{guo35}) and exchanging the order of summation lead to
\begin{equation}\label{guo313}
\begin{aligned}
\phi(z) &=c_{1}+\sum_{k=2}^{n-1} \sum_{j=k+1}^{n}(-1)^{j+1} \frac{2 j-2 k}{2 j-1} C_{n}^{j} z^{2 k-2} \\
&=c_{1}+\sum_{j=3}^{n} \sum_{k=2}^{j-1}(-1)^{j+1}\left(1+\frac{1-2 k}{2 j-1}\right) C_{n}^{j} z^{2 k-2} \\
& \triangleq c_{1}+\tilde{\phi}_{1}(z)+\tilde{\phi}_{2}(z).
\end{aligned}
\end{equation}
One calculates   $\tilde{\phi}_{1}(z)$ as
\begin{equation}\label{guo314}
\begin{aligned}
\tilde{\phi}_{1}(z) &=\sum_{j=3}^{n}(-1)^{j+1} C_{n}^{j}\left(\frac{z^{2}\left(1-z^{2(j-2)}\right)}{1-z^{2}}\right) \\
&=\frac{z^{2}}{1-z^{2}} \sum_{j=3}^{n}(-1)^{j+1} C_{n}^{j}+\frac{z^{-2}}{1-z^{2}} \sum_{j=3}^{n}(-1)^{j} C_{n}^{j} z^{2 j} \\
& \triangleq \tilde{\phi}_{1,1}(z)+\tilde{\phi}_{1,2}(z).
\end{aligned}
\end{equation}
Further computation generates
\begin{equation}\label{guo315}
\begin{aligned}
\tilde{\phi}_{1,1}(z)=-\frac{z^{2}}{1-z^{2}}\left(\sum_{j=0}^{n}(-1)^{j} C_{n}^{j}-\sum_{j=0}^{2}(-1)^{j} C_{n}^{j}\right)
=\frac{z^{2}}{1-z^{2}} \cdot \sum_{j=0}^{2}(-1)^{j} C_{n}^{j}
=\frac{n^{2}-3n+2}{2}\cdot \frac{z^{2}}{1-z^{2}}
\end{aligned}
\end{equation}
and
\begin{equation}\label{guo316}
\begin{aligned}
\tilde{\phi}_{1,2}(z) =\frac{z^{-2}}{1-z^{2}}\left(\sum_{j=0}^{n}\left(-z^{2}\right)^{j} C_{n}^{j}-\sum_{j=0}^{2}\left(-z^{2}\right)^{j} C_{n}^{j}\right)
=\frac{z^{-2}}{1-z^{2}}\left(\left(1-z^{2}\right)^{n}-1+n z^{2}-z^{4} C_{n}^{2}\right) \triangleq \frac{z^{-2}}{1-z^{2}}
\rho(z).
\end{aligned}
\end{equation}
Consequently, combining  Eqs.~(\ref{guo314})-(\ref{guo316}) produces
\begin{equation}\label{guo317}
\tilde{\phi}_{1}(z)=\frac{n^{2}-3n+2}{2} \frac{z^{2}}{1-z^{2}}+\frac{z^{-2}}{1-z^{2}}\rho(z).
\end{equation}
For $\tilde{\phi}_{2}(z)$, direct computation yields
\begin{equation}\label{guo318}
\begin{aligned}
\tilde{\phi}_{2}(z) &=\frac{\mathrm{d}}{\mathrm{d} z}\left(\sum_{j=3}^{n} \sum_{k=2}^{j-1} \frac{(-1)^{j} C_{n}^{j}}{2 j-1} z^{2 k-1}\right) \\
&=\frac{\mathrm{d}}{\mathrm{d}z}\left(\sum_{j=3}^{n} \frac{(-1)^{j} C_{n}^{j}}{2 j-1} \cdot \frac{z^{3}\left(1-z^{2(j-2)}\right)}{1-z^{2}}\right) \\
&=\frac{\mathrm{d}}{\mathrm{d} z}\left(\frac{z^{3}}{1-z^{2}} \cdot \sum_{j=3}^{n} \frac{(-1)^{j} C_{n}^{j}}{2 j-1}-\frac{1}{1-z^{2}} \cdot \sum_{j=3}^{n} \frac{(-1)^{j} C_{n}^{j}}{2 j-1} z^{2 j-1}\right) \\
& \triangleq \frac{\mathrm{d}}{\mathrm{d} z}\left(\tilde{\Phi}_{2,1}(z)+\tilde{\Phi}_{2,2}(z)\right).
\end{aligned}
\end{equation}

Using similar method as (\ref{guo316}), one derives
\begin{equation}\label{guo319}
\begin{aligned}
\tilde{\Phi}_{2,1}(z) &=\frac{z^{3}}{1-z^{2}} \cdot \int_{0}^{1} \sum_{j=3}^{n}(-1)^{j} C_{n}^{j} z^{2 j-2} \mathrm{d} z \\
&=\frac{z^{3}}{1-z^{2}} \cdot \int_{0}^{1} \sum_{j=3}^{n} C_{n}^{j}\left(-z^{2}\right)^{j} \frac{\mathrm{d} z}{z^{2}}
=\frac{z^{3}}{1-z^{2}} \cdot \int_{0}^{1} \frac{\rho(z)}{z^{2}} \mathrm{d} z
\end{aligned}
\end{equation}
and
\begin{equation}\label{guo320}
\tilde{\Phi}_{2,2}(z)=\frac{1}{z^{2}-1} \cdot \int_{0}^{z} \sum_{j=3}^{n}\left(-z^{2}\right)^{j} C_{n}^{j} \frac{\mathrm{d} z}{z^{2}}=\frac{1}{z^{2}-1} \cdot \int_{0}^{z} \frac{\rho(z)}{z^{2}} \mathrm{d} z.
\end{equation}
Gathering  (\ref{guo318})-(\ref{guo320}) leads to
\begin{equation}\label{guo321}
\tilde{\phi}_{2}(z)=\frac{3 z^{2}\d\int_{0}^{1} \frac{\rho(s)}{s^{2}} ds-\frac{\rho(z)}{z^{2}}}{1-z^{2}}
+\frac{2 z\left(z^{3}\d\int_{0}^{1} \frac{\rho(s)}{s^{2}} d s-\int_{0}^{z} \frac{\rho(s)}{s^{2}} d s\right)}{\left(1-z^{2}\right)^{2}}.
\end{equation}
Similar as  (\ref{guo319}), we invoke  (\ref{guo35})
to compute  $c_{1}$ to find
\begin{equation}\label{guo322}
\begin{aligned}
c_{1} &=\frac{1}{2}+\sum_{j=1}^{n}(-1)^{j+1} \frac{C_{n}^{j}}{2}-\sum_{j=1}^{n}(-1)^{j+1} \frac{C_{n}^{j}}{2 j-1} \\
&=1+\sum_{j=0}^{n}(-1)^{j+1} \frac{C_{n}^{j}}{2}+\sum_{j=1}^{n} \frac{(-1)^{j} C_{n}^{j}}{2 j-1}=1-n+\frac{C_{n}^{2}}{3}+\int_{0}^{1} \frac{\rho(s)}{s^{2}} \mathrm{d}s.
\end{aligned}
\end{equation}
Employing the definition of $\omega(z)$ in  (\ref{guo316}), one finds
\begin{equation}\label{guo323}
\begin{aligned}
\int_{0}^{1} \frac{\rho(s)}{s^{2}} \mathrm{d}s &=-\int_{0}^{1} \sum_{k=0}^{n-1}\left(1-z^{2}\right)^{k} \mathrm{d}z
+\int_{0}^{1}\left(n-z^{2} C_{n}^{2}\right) \mathrm{d} z \\
&=-\sum_{k=1}^{n-1}\int_{0}^{1} \left(1-z^{2}\right)^{k} \mathrm{d}z+(n-1)-\frac{C_{n}^{2}}{3} \triangleq-B+(n-1)-\frac{C_{n}^{2}}{3}.
\end{aligned}
\end{equation}

It follows from (\ref{guo313}), (\ref{guo317}) and (\ref{guo321})-(\ref{guo323}) that one has
\begin{equation*}
\begin{aligned}
\phi(z)=&-B+\frac{(2 (n-1)-3 B) z^{2}}{1-z^{2}}+\frac{2z^{4}\left(-B+(n-1)-\frac{1}{3}C_{n}^{2}\right)}
{\left(1-z^{2}\right)^{2}}+\frac{2 z\d\int_{0}^{z}\left(\sum_{k=0}^{n-1}\left(1-s^{2}\right)^{k}-n+s^{2} C_{n}^{2}\right) \mathrm{d} s}{\left(1-z^{2}\right)^{2}} \\
=& \frac{-B-B z^{2}+2 z \d\int_{0}^{z} \sum_{k=1}^{n-1}\left(1-s^{2}\right)^{k}\mathrm{d}s}
{\left(1-z^{2}\right)^{2}} \\
\leq & \frac{-B-B z^{2}}{\left(1-z^{2}\right)^{2}}+\frac{2 z   B}
{\left(1-z^{2}\right)^{2}}=-\frac{B}{(1+z)^{2}} \leq 0
\end{aligned}
\end{equation*}
implying  (\ref{guo312}) for $z \in[0,1)$.
Then using  (\ref{guo35}) and (\ref{guo44}), after  a similar discussion,
one finds  that
$(\ref{guo39})$ is equivalent to $\sum_{k=1}^{n-1} c_{2 k-1} z^{2 k-2}(z-1)^{2} \leq 0$,
 which is obviously true by (\ref{guo312}).
Therefore, we complete the proof of  (\ref{guo38}).
Gathering  (\ref{guo34}), (\ref{guo36}) and (\ref{guo38}) generates
\begin{equation*}
\begin{aligned}
&F(u)-\frac{2(2-c_{1})}{2n+1} M^{2 n+1}(t)
=\int_{\mathbb{R}} h(t, x) g^{2}(t, x) \mathrm{d} x \\
&\qquad \leq \frac{2-c_{1}}{2} M^{2 n-1}(t) \int_{\mathbb{R}} g^{2}(t, x) \mathrm{d} x
=\frac{2-c_{1}}{2}\left(E(u)-2 M^{2}(t)\right)M^{2 n-1}(t)
\end{aligned}
\end{equation*}
indicating  (\ref{guo33}). This completes the proof of Proposition \ref{lemGuo32}.
\end{proof}

\begin{proof}[\textit{\textbf{Proof of Proposition \ref{lemGuo33}}}]
Since  $0<\varepsilon<(\gamma-2 \sqrt{2}) a,\, \, \gamma>2\sqrt{2}$ and $\left\|\varphi_{c}\right\|_{H^{1}}^2=2a^2$ given by Eq.~(\ref{Epeakon}), one easily deduces
\begin{equation}\label{guo324}
\begin{aligned}
\left|E(u)-E\left(\varphi_{c}\right)\right| &=\left|\left(\|u\|_{H^{1}}-\left\|\varphi_{c}\right\|_{H^{1}}\right)\left(\|u\|_{H^{1}}+\left\|\varphi_{c}\right\|_{H^{1}}\right)\right| \\
& \leq\left\|u-\varphi_{c}\right\|_{H^{1}}\left(\left\|u-\varphi_{c}\right\|_{H^{1}}+2\left\|\varphi_{c}\right\|_{H^{1}}\right) \leq \varepsilon(\varepsilon+2 \sqrt{2} a) < a \gamma \varepsilon.
\end{aligned}
\end{equation}
Moreover, we find
\begin{equation}\label{guo325}
\begin{aligned}
&\left|F(u)-F\left(\varphi_{c}\right)\right| \\
&= \bigg|\int_{\mathbb{R}}\left(u^{2 n+1}+\sum_{k=1}^{n} \frac{(-1)^{k+1}}{2 k-1} C_{n}^{k} u^{2 n-2 k+1} u_{x}^{2 k}\right) \mathrm{d} x
-\int_{\mathbb{R}}\left(\varphi_{c}^{2 n+1}+\sum_{k=1}^{n} \frac{(-1)^{k+1}}{2 k-1} C_{n}^{k} \varphi_{c}^{2 n-2 k+1}\left(\partial_{x} \varphi_{c}\right)^{2 k} \right) \mathrm{d} x \bigg| \\
&\leq \bigg|\int_{\mathbb{R}}\left(u^{2 n+1}+n u^{2 n-1} u_{x}^{2}-\varphi_{c}^{2 n+1}-n \varphi_{c}^{2 n-1}\left(\partial_{x} \varphi_{c}\right)^{2}\right) \mathrm{d} x \bigg| \\
& \qquad+\sum_{k=2}^{n} \frac{C_{n}^{k}}{2 k-1}\bigg|\int_{\mathbb{R}}\left(u^{2 n-2 k+1} u_{x}^{2 k}-\varphi_{c}^{2 n-2 k+1}\left(\partial_{x} \varphi_{c}\right)^{2 k}\right) \mathrm{d} x\bigg|\\
&\leq \bigg|\int_{\mathbb{R}}\left(u^{2 n-1}-\varphi_{c}^{2 n-1}\right)\left(u^{2}
+n u_{x}^{2}\right) \mathrm{d} x\bigg|
+\left|\int_{\mathbb{R}} \varphi_{c}^{2 n-1}\left(\left(u^{2}-\varphi_{c}^{2}\right)+n\left(u_{x}^{2}-\left(\partial_{x} \varphi_{c}\right)^{2}\right)\right) \mathrm{d} x\right| \\
&\qquad+\sum_{k=2}^{n} \frac{C_{n}^{k}}{2 k-1} \bigg| \int_{\mathbb{R}} u^{2 n-2 k+1}
\left(u_{x}^{2 k}-\left(\partial_{x} \varphi_{c}\right)^{2 k}\right) \mathrm{d} x\bigg|+\sum_{k=2}^{n} \frac{C_{n}^{k}}{2 k-1}\bigg| \int_{\mathbb{R}}\left(\partial_{x} \varphi_{c}\right)^{2 k}\left(u^{2 n-2 k+1}-\varphi_{c}^{2 n-2 k+1}\right) \mathrm{d} x \bigg| \\
&\triangleq A_{1}+A_{2}+A_{3}+A_{4}.
\end{aligned}
\end{equation}

We will estimate $A_{1}$-$A_{4}$ term by term.
We employ Eq.~(\ref{guo324}) to control the term $A_{1}$ as
\begin{equation*}
\begin{aligned}
A_{1} \leq & \, n \int_{\mathbb{R}}\left|u-\varphi_{c}\right| \cdot\left|u^{2 n-2}+u^{2 n-3} \varphi_{c}+\cdots+u \varphi_{c}^{2 n-2}+\varphi_{c}^{2 n-2}\right| \cdot\left(u^{2}+u_{x}^{2}\right) \mathrm{d} x\\
\leq &\, n E(u)\left\|u-\varphi_{c}\right\|_{L^{\infty}}\left(\|u\|_{L^{\infty}}^{2 n-2}
+\|u\|_{L^{\infty}}^{2 n-3}\left\|\varphi_{c}\right\|_{L^{\infty}}
+\cdots+\|u\|_{L^{\infty}}\left\|\varphi_{c}\right\|_{L^{\infty}}^{2 n-3}
+\left\|\varphi_{c}\right\|_{L^{\infty}}^{2 n-2}\right) \\
\leq &\, \frac{n}{2^{n-1}}\left(E\left(\varphi_{c}\right)+a \gamma \varepsilon\right)\left\|u-\varphi_{c}\right\|_{H^{1}}\left(\|u\|_{H^{1}}^{2 n-2}+\|u\|_{H^{1}}^{2 n-3}\left\|\varphi_{c}\right\|_{H^{1}}+\cdots
+\|u\|_{H^{1}}\left\|\varphi_{c}\right\|_{H^{1}}^{2 n-3}
+\left\|\varphi_{c}\right\|_{H^{1}}^{2 n-2}\right) \\
\leq &\, \frac{n}{2^{n-1}}\left(E\left(\varphi_{c}\right)+ a \gamma \varepsilon\right)\left(\left\|u-\varphi_{c}\right\|_{H^{1}}
+2\left\|\varphi_{c}\right\|_{H^{1}}\right)^{2 n-2} \left\|u-\varphi_{c}\right\|_{H^{1}} \\
\leq &\, \frac{n}{2^{n-1}}\varepsilon\left(2 a^{2}+ a \gamma\varepsilon\right)(\varepsilon+2 \sqrt{2} a)^{2 n-2},
\end{aligned}
\end{equation*}
where we have also used the following inequality
\begin{eqnarray}\label{guo46}
\sup _{x \in \mathbb{R}}|v(x)| \leq \frac{\sqrt{E(v)}}{\sqrt{2}} \leq \frac{\|v\|_{H^{1}}}{\sqrt{2}},\quad \text{for}\quad v\in H^{1}(\mathbb{R})
\end{eqnarray}
coming from (\ref{guo34}).

A similar argument applied  for  $A_{2}$ lead to
\begin{equation*}
\begin{aligned}
A_{2} \leq &\, \left\|\varphi_{c}\right\|_{L^{\infty}}^{2 n-1} \bigg| \int_{\mathbb{R}}\left(\left(u-\varphi_{c}\right)^{2}+n\left(u_{x}-\partial_{x} \varphi_{c}\right)^{2}+2 \varphi_{c}\left(u-\varphi_{c}\right)
+2n \partial_{x} \varphi_{c}\left(u_{x}-\partial_{x} \varphi_{c}\right)\right) \mathrm{d} x \bigg| \\
\leq &\,\left(\frac{\left\|\varphi_{c}\right\|_{H^{1}}}{\sqrt{2}}\right)^{2 n-1}\left((n+1)\left\|u-\varphi_{c}\right\|_{H^{1}}^{2}+2(n+1)\left\|\varphi_{c}\right\|_{H^{1}}\left\|u-\varphi_{c}\right\|_{H^{1}}\right) \\
\leq &\, (n+1) a^{2 n-1}\varepsilon(\varepsilon+2 \sqrt{2} a).
\end{aligned}
\end{equation*}
The term $A_{3}$ can be handled by
\begin{equation*}
\begin{aligned}
A_{3} \leq &\, \sum_{k=2}^{n} \frac{C_{n}^{k}}{2 k-1}\|u\|_{L^{\infty}}^{2 n-2 k+1}\left(\int_{\mathbb{R}} ( u_{x}^{k} + ( \partial_{x} \varphi_{c} )^{k})^{2} \left(u_{x}^{k-1}+u_{x}^{k-2} \partial_{x} \varphi_{c}+\cdots+u_{x}\left(\partial_{x} \varphi_{c}\right)^{k-2}\right.\right.\\
&\,\left.\left.+\left(\partial_{x} \varphi_{c}\right)^{k-1}\right)^{2} \mathrm{d} x\right)^{\frac{1}{2}}\left(\int_{\mathbb{R}}\left(u_{x}-\partial_{x} \varphi_{c}\right)^{2} \mathrm{d} x\right)^{\frac{1}{2}} \\
=&\, \sum_{k=2}^{n} \frac{C_{n}^{k}}{2 k-1}\|u\|_{L^{\infty}}^{2 n-2 k+1}\left(\int _ { \mathbb { R } } \left(u_{x}^{4 k-2}+2 u_{x}^{4 k-3} \partial_{x} \varphi_{c}+\cdots+(2 k-1) u_{x}^{2 k}\left(\partial_{x} \varphi_{c}\right)^{2 k-2}\right.\right.\\
&\,+2 k u_{x}^{2 k-1}\left(\partial_{x} \varphi_{c}\right)^{2 k-1}+(2 k-1) u_{x}^{2 k-2}\left(\partial_{x} \varphi_{c}\right)^{2 k}+\cdots+2 u_{x}\left(\partial_{x} \varphi_{c}\right)^{4 k-3} \\
&\,\left.\left.+\left(\partial_{x} \varphi_{c}\right)^{4 k-2}\right) \mathrm{d} x\right)^{\frac{1}{2}}\left(\int_{\mathbb{R}}\left(u_{x}-\partial_{x} \varphi_{c}\right)^{2} \mathrm{d} x\right)^{\frac{1}{2}}\\
\leq & \,\sum_{k=2}^{n} \frac{C_{n}^{k}}{2 k-1}\|u\|_{L^{\infty}}^{2 n-2 k+1} \cdot \sqrt{2} k\left(\int_{\mathbb{R}} u_{x}^{4 k-2} \mathrm{d} x+\int_{\mathbb{R}}\left(\partial_{x} \varphi_{c}\right)^{4 k-2} \mathrm{d} x\right)^{\frac{1}{2}} \cdot\left\|u-\varphi_{c}\right\|_{H^{1}} \\
\lesssim &\, \|u\|_{H^{1}}^{2 n-2 k+1}\left(\|u\|_{L^{2}}^{k}\left\|u_{x x}\right\|_{L^{2}}^{3 k-2}+\left\|\partial_{x} \varphi_{c}\right\|_{L^{4 k-2}}^{4 k-2}\right)^{\frac{1}{2}} \cdot\left\|u-\varphi_{c}\right\|_{H^{1}} \\
\lesssim &\, G\left(\|u\|_{H^{s}},n, c\right) \varepsilon,
\end{aligned}
\end{equation*}
where we have used the Gagliardo-Nirenberg inequality
$\left\|u_{x}\right\|_{L^{4 k-2}}^{4 k-2}
 \leq C\|u\|_{L^{2}}^{k}\left\|u_{x x}\right\|_{L^{2}}^{3 k-2}$
and the fact $\left\|\partial_{x} \varphi_{c}\right\|_{L^{4 k-2}}^{4 k-2}
=\frac{a^{4 k-2}}{2 k-1}$.

For the term $A_{4}$,  the H\"{o}lder inequality produces
\begin{equation*}
\begin{aligned}
A_{4} \leq & \sum_{k=2}^{n} \frac{C_{n}^{k}}{2 k-1}
\left(\|u\|_{L^{\infty}}^{2n-2k}+\|u\|_{L^{\infty}}^{2n-2k-1}
\left\|\varphi_{c}\right\|_{L^{\infty}}+\cdots
+\|u\|_{L^{\infty}}\left\|\varphi_{c}\right\|_{L^{\infty}}^{2n-2k-1}
+\left\|\varphi_{c}\right\|_{L^{\infty}}^{2n-2k}\right)\\
&\quad\quad\times\left(\int_{\mathbb{R}}\left(\partial_{x} \varphi_{c}\right)^{4 k} \mathrm{d} x\right)^{\frac{1}{2}}\left(\int_{\mathbb{R}}\left(u-\varphi_{c}\right)^{2} \mathrm{d} x\right)^{\frac{1}{2}} \\
\leq & \sum_{k=2}^{n} \frac{C_{n}^{k}}{2 k-1}\left(\left\|u-\varphi_{c}\right\|_{H^{1}}+2\left\|\varphi_{c}\right\|_{H^{1}}\right)^{2 n-2 k}\left\|\partial_{x} \varphi_{c}\right\|_{L^{4 k}}^{2 k} \cdot\left\|u-\varphi_{c}\right\|_{H^{1}} \\
\leq & \sum_{k=2}^{n} \frac{C_{n}^{k}a^{4 k}}{2k(2k-1)}\varepsilon(2 \sqrt{2} a+\varepsilon)^{2(n-k)}.
\end{aligned}
\end{equation*}
Substituting the above estimates of $K_{1}$-$K_{4}$ into (\ref{guo325}) yields
\begin{eqnarray*}
(\ref{guo325})\lesssim G\left(\|u\|_{H^{s}}, n, c\right) \cdot \varepsilon.
\end{eqnarray*}
 This completes the proof of Proposition \ref{lemGuo33}.
\end{proof}

\begin{proof}[\textit{\textbf{Proof of Proposition \ref{lemGuo34}}}]
We first derive from (\ref{guo33}) that
\begin{eqnarray}\label{guo326}
(2n-1) M^{2 n+1}-\frac{2n+1}{2} M^{2 n-1} E(u)+\frac{2n+1}{2-c_{1}} F(u) \leq 0.
\end{eqnarray}
Let
 $Q(y)=(2n-1) y^{2 n+1}-\frac{2n+1}{2} y^{2 n-1} E(u)+\frac{2n+1}{2-c_{1}} F(u).$
For the case of $E(u)=E\left(\varphi_{c}\right)=2 a^{2}$ and $F(u)=F\left(\varphi_{c}\right)=\frac{2}{2n+1} a^{2 n+1}\left(1+\sum_{k=1}^{n} \frac{(-1)^{k+1}}{2 k-1} C_{n}^{k}\right)$, $Q(y)$ reduces to
\begin{equation}\label{guo327}
\begin{aligned}
\hat Q(y) &=(2n-1) y^{2 n+1}
-\frac{2n+1}{2} E\left(\varphi_{c}\right) y^{2 n-1}
+\frac{2n+1}{2-c_{1}} F\left(\varphi_{c}\right) \\
&=(2n-1) y^{2 n+1}-(2n+1) a^{2} y^{2 n-1}+\frac{2}{2-c_{1}} \cdot a^{2 n+1}\left(1+\sum_{k=1}^{n} \frac{(-1)^{k+1}}{2 k-1} C_{n}^{k}\right)\\
&=(2n-1) y^{2 n+1}-(2n+1) a^{2} y^{2 n-1}+2a^{2 n+1} \\
&=(y-a)^{2}\left((2n-1) y^{2 n-1}+2\sum_{k=1}^{2 n-2}(2 n-k) a^{k} y^{2 n-1-k}
+2a^{2 n-1}\right),
\end{aligned}
\end{equation}
where the relation
\begin{equation}\label{guo329}
\begin{aligned}
2-c_{1}=1+\sum_{k=1}^{n} \frac{(-1)^{k+1}}{2 k-1} C_{n}^{k}.
\end{aligned}
\end{equation}
is used since $\sum_{k=1}^{n}(-1)^{k+1} C_{n}^{k}=1$.

Then, it follows from (\ref{guo326}) and (\ref{guo327}) that one has
\begin{equation*}
\hat Q(M) \leq \hat Q(M)-Q(M)
=\frac{2n+1}{2} M^{2 n-1}
\left(E(u)-E\left(\varphi_{c}\right)\right)
-\frac{2n+1}{2-c_{1}}\left(F(u)-F\left(\varphi_{c}\right)\right),
\end{equation*}
which along with  (\ref{guo327})  generates
\begin{equation}\label{guo331}
2a^{2 n-1}(M-a)^{2} \leq \frac{2n+1}{2} M^{2 n-1}
\left(E(u)-E\left(\varphi_{c}\right)\right)
-\frac{2n+1}{2-c_{1}}\left(F(u)-F\left(\varphi_{c}\right)\right).
\end{equation}
 Thanks to  (\ref{guo34}) and the assumption of this lemma, one derives
  for $0<\varepsilon<(\gamma-2 \sqrt{2}) a\,\, (\gamma>2\sqrt{2})$
  that
\begin{equation}\label{guo332}
0<M^{2} \leq \frac{E(u)}{2} \leq \frac{2 a^{2}+ a \gamma \varepsilon}{2}<\frac{(\gamma-\sqrt{2})a^{2}}{2}.
\end{equation}
Gathering  (\ref{guo331}) and  (\ref{guo332}) gives rise to
\begin{equation*}
\sqrt{2}a^{n-1/2}|M-a| \lesssim \sqrt{\frac{(2n+1) \gamma(\gamma-\sqrt{2})^{n}}{2^{n}\sqrt{2(\gamma-\sqrt{2})}} a^{2 n} \varepsilon+\frac{2n+1}{2} G\left(n, c,\|u\|_{H^{s}}\right) \varepsilon}.
\end{equation*}
Combing  this inequality and the relation
   $c=\sum_{k=0}^{n}\frac{(-1)^{k}}{2k+1}\frac{2n+1}{2n}C_{n}^{k}a^{2n-1}$,
   one concludes that there exists  a constant, still expressed via $G\left(n, c,\|u\|_{H^{s}}\right)$, such that
\begin{equation*}
|M-a| \lesssim \sqrt{G(n, c,\|u\|_{H^{s}})\varepsilon}.
\end{equation*}
This completes the proof of Proposition \ref{lemGuo34}.
\end{proof}

We finally prove Theorem \ref{thmGuo31}.

\begin{proof}[\textit{\textbf{Proof of Theorem  \ref{thmGuo31}}}]
Let $u \in C\left([0, T) ; H^{s}(\mathbb{R})\right)(s>\frac{5}{2})$ satisfy Eq.~(\ref{gCH}) with the  initial condition $u_{0} \in H^{s}(\mathbb{R})$. Moreover one can know
\begin{equation}\label{guo333}
E(u(t, \cdot))=E\left(u_{0}\right), \qquad  F(u(t, \cdot))=F\left(u_{0}\right), \quad \forall t \in[0, T).
\end{equation}
Since $\left\|u(0, \cdot)-\varphi_{c}\right\|_{H^{s}}<\varepsilon$, with $0<\varepsilon<(\gamma-2 \sqrt{2}) a$, and $0 \neq\left(1-\partial_{x}^{2}\right) u_{0} \geq 0$, in view of (\ref{guo333}) and Proposition \ref{lemGuo33}, the hypotheses of Proposition \ref{lemGuo34} are satisfied for $u(t, \cdot)$ with a chosen positive constant $G\left(n, c,\left\|u_{0}\right\|_{H^{s}}\right)$ depending only on the wave speed $c$, $n\in\mathbb{N}^+$ and $\left\|u_{0}\right\|_{H^{s}}$. Accordingly
\begin{equation}\label{guo334}
|u(t, \eta(t))-a| \lesssim \sqrt{G\left(n, c,\left\|u_{0}\right\|_{H^{s}}\right) \varepsilon}, \quad \forall t \in[0, T).
\end{equation}
where $x=\eta(t)\in\mathbb{R}$ stands for the maximum point of function $u(t,x)$. Consequently, we conclude from  (\ref{guo334}) and Proposition \ref{lemGuo33} that for $t \in[0, T)$,
$$
\left\|u-\varphi_{c}(\cdot-\eta(t))\right\|_{H^{1}} \leq \sqrt{\left|E\left(u_{0}\right)-E\left(\varphi_{c}\right)\right|+4 a|u(t, \eta(t))-a|} \lesssim \sqrt{a \gamma \varepsilon+4 a \sqrt{G\left(c,\left\|u_{0}\right\|_{H^{s}}\right) \varepsilon}}.
$$
We thus complete  the proof of Theorem \ref{thmGuo31}.
\end{proof}

\v\v\textbf{Acknowledgement}

This work was partially supported by the NSFC of China under Grants No. 11925108.

\begin{appendix}

\vspace{0.3in}
\noindent {\bf Appendix A.\, Some basics properties of Eq.~(\ref{gCH}) with the initial condition $u(0, x)=u_0(x)$ } \v

\setcounter{equation}{0}
\renewcommand\theequation{A.\arabic{equation}}

\newtheorem{lemma}{Lemma}
\renewcommand\thelemma{A.\arabic{lemma}}

\begin{lemma}\cite{QuFu2020JDDE}\label{lemGuo21}
 Let $u_{0}(x)=u(0, x) \in H^{s}(\mathbb{R})$ with $s>5 / 2$. Then there exists a time $T>0$ such that the Cauchy problem of Eq.~(\ref{gCH}) with initial data $u_{0}$
 has a unique strong solution $u(t, x) \in C\left([0, T) ; H^{s}(\mathbb{R})\right) \cap$ $C^{1}\left([0, T) ; H^{s-1}(\mathbb{R})\right)$ and the map $u_{0} \mapsto u$ is continuous from a neighborhood of $u_{0}$ in $H^{s}(\mathbb{R})$ into $C\left([0, T) ; H^{s}(\mathbb{R})\right) \cap C^{1}\left([0, T) ; H^{s-1}(\mathbb{R})\right)$.
\end{lemma}
For the characteristic equation
\begin{equation}\label{guo21}
\left\{\begin{array}{l}
\dfrac{d r(t, x)}{d t}=\left(u\left(u^{2}-u_{x}^{2}\right)^{n-1}\right)(t, r(t, x)), x \in \mathbb{R}, \quad t \in[0, T) \\
r(0, x)=x, x \in \mathbb{R},
\end{array}\right.
\end{equation}
there exists the following lemma: 

\begin{lemma}\cite{QuFu2020JDDE}\label{propGuo1}
 Suppose $u_{0} \in H^{s}(\mathbb{R})$ with $s>5 / 2$, and let $T>0$ be the maximal existence time of the strong solution $u(t, x) \in C\left([0, T), H^{s}(\mathbb{R})\right) \cap C^{1}\left([0, T), H^{s-1}(\mathbb{R})\right)$ to the Cauchy problem (\ref{gCH}). Then (\ref{guo22}) has a unique solution $r \in C^{1}([0, T) \times \mathbb{R}, \mathbb{R})$ such that $r(t, \cdot)$ is an increasing diffeomorphism over $\mathbb{R}$ with
\begin{eqnarray}\label{guo22}
r_{x}(t, x)=\exp \left(\int_{0}^{t}\left(u_{x}\left(u^{2}-u_{x}^{2}\right)^{n-1}+2 (n-1) u u_{x}\left(u^{2}-u_{x}^{2}\right)^{n-2} m\right)(s, r(s, x)) \mathrm{d} s\right)>0
\end{eqnarray}
for all $(t, x) \in[0, T) \times \mathbb{R}$.

Furthermore, the momentum density $m=u-u_{x x}$ satisfies
\begin{eqnarray}\label{guo22}
m(t, r(t, x)) r_{x}(t, x)=m_{0}(x) \exp \left[-\int_{0}^{t}\left(u_{x}\left(u^{2}-u_{x}^{2}\right)^{n-1}\right)(s, r(s, x)) \mathrm{d} s\right],
\end{eqnarray}
for all $(t, x) \in[0, T) \times \mathbb{R}$,
which implies that the sgn and zeros of $m$ are preserved under the flow.
\end{lemma}
Based on the above Lemma, there was the following property: 
\begin{lemma}\cite{QuFu2020JDDE}\label{lemGuo22}
Assume $u_{0}(x) \in H^{s}(\mathbb{R}), s>5 / 2$. If $m_{0}(x)=\left(1-\partial_{x}^{2}\right) u_{0}(x)$ does not change sgn, then $y(t, x)$ will not change sgn for all $t \in[0, T) .$ It follows that if $m_{0} \geq 0$, then the corresponding solution $u(t, x)$ of equation (\ref{gCH}) is positive for $(t, x) \in[0, T) \times \mathbb{R}$. Furthermore, if $m_{0} \geq 0$, then the corresponding solution u $(t, x)$ of equation (\ref{gCH}) satisfies
$$
\left(1 \pm \partial_{x}\right) u(t, x) \geq 0, \quad \text { for } \quad \forall(t, x) \in[0, T) \times \mathbb{R}.
$$
\end{lemma}

\end{appendix}

\end{document}